# Method for Undershoot-Less Control of Non-Minimum Phase Plants Based on Partial Cancellation of the Non-Minimum Phase Zero: Application to Flexible-Link Robots

F. Merrikh-Bayat and F. Bayat
Department of Electrical and Computer Engineering
University of Zanjan
Zanjan, Iran
Email: f.bayat@znu.ac.ir , bayat.farhad@znu.ac.ir

*Abstract*—As a well understood classical fact, non- minimum phase zeros of the process located in a feedback connection cannot be cancelled by the corresponding poles of controller since such a cancellation leads to internal instability. This impossibility of cancellation is the source of many limitations in dealing with the feedback control of non-minimum phase processes. The aim of this paper is to study the possibility and usefulness of partial (fractional-order) cancellation of such zeros for undershoot-less control of non-minimum phase processes. In this method first the non-minimum phase zero of the process is cancelled to an arbitrary degree by the proposed pre-compensator and then a classical controller is designed to control the series connection of these two systems. Since plants with multiple non-minimum phase zeros and oscillatory poles are very common in the problems related to robotics, the proposed method is applied to these systems to confirm its effectiveness.

I. INTRODUCTION

It is well understood that non-minimum phase processes constitute a challenging research area in the field of control engineering. Non-minimum phase zeros appear unavoidably in some important industrial processes such as steam generators [1], aircrafts [2], [3], flexible-link manipulators [4], continuous stirred tank reactors (CSTRs) [5], electronic circuits [6], and so on. As a very well known classical fact, non-minimum phase zeros of the process put some limitations on the performance of the corresponding feedback system [7]-[10]. More precisely, these limitations can be concluded, e.g., from the classical root-locus method [11], asymptotic LQG theory [9], waterbed effect phenomena [12], and the LTR problem [13]. In the field of linear time-invariant (LTI) systems, the source of all of the above-mentioned limitations is that the non-minimum phase zero of the given process cannot be cancelled by unstable pole of the controller since such a cancellation leads to internal instability [14].

During the past decades various methods have been developed by researchers for the control of processes with non-minimum phase zeros (see, for example, [15]-[17] and the references therein for more information on this subject). Among others, according to the simplicity and high achievement of the feedback control strategy in dealing with most of the real-world industrial problems, it is strictly preferred to develop more effective methods to the control of non-minimum phase processes by using this technique. However, as mentioned before, impossibility of unstable pole-zero cancellation is the main limitation of this method, which is to be partly removed in this paper.

An author of this paper already showed [18] that although unstable pole-zero cancellation is impractical in LTI feedback systems and leads to internal instability, the *partial* (or, *fractional-order*) unstable pole-zero cancellation is possible and can be very effective. In fact, it is proved in [18] that any non-minimum phase zero (unstable pole) of the given process can partly be cancelled by a pole (zero) of the controller without resulting in an internally-unstable feedback system. Interesting observation is that this cancellation can also increase the phase and gain margin of the closed-loop system, and consequently, partly remove some of the classical limitations caused by non-minimum phase zeros. Note that the method proposed in [18] can be used to cancel any non-minimum phase zero or unstable pole of a process to an arbitrary degree.

The aim of this paper is to study the control of certain class of robot arms by combining the proposed method for cancellation of non-minimum phase zeros of the process and the classical PID control. A relatively similar approach, which studies the integral performance indices of a feedback system (in which a PI controller is applied in series with a fractional-order pole-zero canceller to control a second order process) is presented in [19]. Here it is worth to mention that PID controllers commonly do not lead to satisfactory results when the process is non-minimum phase, has poles with a very low damping ratio, or exhibit large dead times [20]. Hence, from the practical point of view it is very important to develop effective methods to remove these limitations. Since transfer functions with multiple non-minimum phase zeros and oscillatory poles frequently appear in dealing with flexible arm robots, the studies of this paper are mainly focused on these systems. However, the proposed ideas are very general and can be applied to any other non-minimum phase process as well.

The rest of this paper is organized as follows. The proposed method for the control of non-minimum phase processes is presented in Section II. Illustrative examples, which are adopted from flexible-link robots, are studied in Section III, and finally Section IV concludes the paper.

## II. MAIN RESULTS

Fig. 1 shows the proposed feedback strategy to control a non-minimum phase process with transfer function $G(s)$ ($r(t)$, $d(t)$ and $n(t)$ stand for the command, disturbance and noise, respectively). As it is observed, in this method first we partially cancel the non-minimum phase zero (or, if necessary, the unstable pole) of $G(s)$ by putting a pre-compensator with transfer function $C_1(s)$ in series with it (see the discussion below). In fact, the role of $C_1(s)$ in Fig. 1 is to remove some of the limitations caused by non-minimum phase zeros of the process by partially removing them. It means that applying $C_1(s)$ will make the control problem *easier* to solve by increasing the phase and gain margin [18].

As it will be shown in the following, $C_1(s)$ is a rational function in non-integer (fractional) powers of $s$. Hence, $P(s) \triangleq C_1(s)G(s)$ in Fig. 1 is a rational function in non-integer powers of $s$ as well. $C_2(s)$ in this figure is used to control the system with transfer function $P(s)$. Note that since $P(s)$ contains fractional powers of $s$, $C_2(s)$ may be designed using either classical design algorithms or the methods specially developed for the control of fractional-order processes (see, for example, [21]-[24] and the references therin for more information on the latter case). For the sake of simplicity we will use the first approach in this paper. In the following, we briefly review the main properties of the fractional-order pole zero canceller, $C_1(s)$, without presenting the proofs. More details can be found in [18].

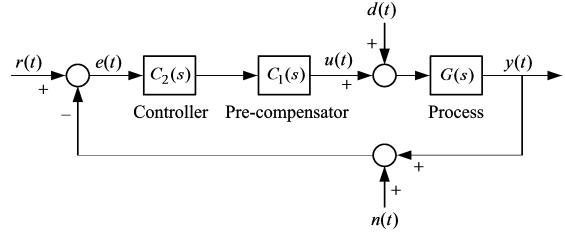

Fig. 1 The general form of the proposed feedback system with pre-compensator (fractional-order pole-zero canceller)

Suppose that $G(s)$ has a positive real zero of order one at $s = z$, that is $G(z) = 0$ and $G'(z) \neq 0$ where $z$ is a positive real number. Such a transfer function can be decomposed as the following:

$$G(s) = \left(1 - \frac{s}{z}\right)\tilde{G}(s). \qquad (1)$$

Clearly, the feedback system shown in Fig. 1 is internally unstable if a pole of $C_1(s)$ (or $C_2(s)$) cancels the non-minimum phase zero of $G(s)$. The following method can be used for partial cancellation of the non-minimum phase zero of $G(s)$ without leading to internal instability. In order to determine the transfer function of the fractional-order pole-zero canceller, $C_1(s)$, first note that the term $1 - s/z$ in (1) can be expanded using fractional powers of $s$ in infinite many different ways as the following:

$$1 - \frac{s}{z} = 1 - \left(\frac{s}{z}\right)^{v/v} = \left[1 - \left(\frac{s}{z}\right)^{1/v}\right]\sum_{k=1}^{v}\left(\frac{s}{z}\right)^{(k-1)/v}, \qquad (2)$$

Where $v$ theoretically can be considered equal to any positive integer. Assuming

$$Q_{z,v}(s) \triangleq \sum_{k=1}^{v}(s/z)^{(k-1)/v}, \qquad (3)$$

Transfer function of the fractional-order unstable pole-zero canceller in Fig. 1 can be defined as the following:

$$C_1(s) = \frac{1}{Q_{z,v}(s)} = \frac{1}{\sum_{k=1}^{v}\left(\frac{s}{z}\right)^{(k-1)/v}}. \qquad (4)$$

(See [25] for time-domain interpretation of fractional powers of $s$ and some real-world examples.) Note that by using the above definition for $C_1(s)$, numerator of the series connection of $C_1(s)$ and $G(s)$ (denoted as $P(s)$) will contain the term $1-(s/z)^{1/v}$ (instead of the term $1-s/z$ in the numerator of $G(s)$), that is

$$P(s) = C_1(s)G(s) = \left[1-\left(\frac{s}{z}\right)^{1/v}\right]\tilde{G}(s). \quad (5)$$

It is proved in [18] that choosing $C_1(s)$ as given in (4), and consequently, changing the non-minimum phase term from $1-s/z$ to $1-(s/z)^{1/v}$ can highly increase the phase and gain margin and partly remove the limitations put on the performance of the feedback system by non-minimum phase zero of the process (of course, without leading to internal instability).

The only unknown parameter of the pre-compensator in Fig. 1 is $v$, which is larger than unity and should be determined by a simple trial and error. Theoretically, the non-minimum phase zero of $G(s)$ can completely be cancelled by tending $v$ to infinity, which is obtained at the cost of using a more complicated setup. However, the problem with applying larger values of $v$ is that it decreases the bandwidth of the open-loop system, and consequently, increases the use of control effort. In practice, in order to design the feedback system first we assign a value to $v$ and then design the controller $C_2(s)$ using a desired method, and next simulate the system. If the responses were satisfactory, the job is done. Else, we should increase the value of $v$ and repeat the procedure.

In general, the controller $C_2(s)$ in Fig. 1 can be designed using any classical controller design algorithm. In this paper $C_2(s)$ is considered as a PID and the effect of the fractional-order pole-zero canceller given in (4) on time-domain responses is studied. Another important alternative for the PID controller to be used in this system is the so-called fractional-order PID (FOPID) or PI$^\lambda$D$^\mu$ controller [24], which is defined as the following:

$$C_2(s) = k_p + \frac{k_i}{s^\lambda} + k_d s^\mu, \quad \lambda, \mu \in \mathbb{R}^+, k_p, k_i, k_d \in \mathbb{R}. \quad (6)$$

Note that unlike classical PID controllers, the FOPID controller given in (6) has five parameters to tune, which makes it a powerful tool to deal with complicated control problems.

According to the above discussions, the feedback system shown in Fig. 2 can be used to control a non-minimum phase process with transfer function $G(s)$. If $G(s)$ has more than one non-minimum phase zero, say at $z_1, \ldots, z_M$, the transfer function of $C_1(s)$ in Fig. 1 should be considered as $1/\prod_{i=1}^{M} Q_{z_i, v_i}(s)$ [18] (see Example 2 of Section III for more details). Note that in this case non-minimum phase zeros can be cancelled to dissimilar degrees, i.e., it is not necessary to subject all of the non-minimum phase zeros of the process to the same amount of cancellation. This technique can also be used for partial cancellation of unstable poles of $G(s)$ [18], which is not discussed in this paper.

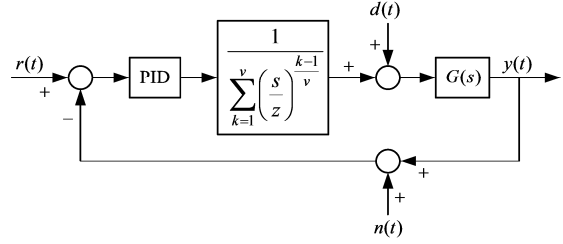

Fig. 2 The feedback system shown in Fig. 1 with a special fractional-order pre-compensator and a PID controller

The last point in relation to the proposed fractional-order pole-zero canceller is about its realization. It general two different methods can be used for this purpose. First, we can approximate the transfer function of $C_1(s)$ with an integer-order transfer function in the frequency range of interest and then realize it using classical methods. The second possible approach is to use the methods available for direct realization of fractional-order systems. See [26]-[29] for more information on the latter case.

### III. ILLUSTRATIVE EXAMPLES

Two illustrative examples are studied in this section to verify the theoretical results of previous section. The processes under consideration in both of these examples are adopted from the problems related to robotics. Since the transfer functions appear in robotics are often non-minimum phase and commonly have oscillatory poles and zeros, they are best suited to the proposed method.

All of the following simulations are performed by taking the numerical inverse Laplace transform from the corresponding transfer functions. More precisely, in each case the unit step response of the feedback system is calculated by taking the numerical inverse Laplace transform from the closed-loop transfer function multiplied by $1/s$. This method is based on the formula proposed in [30] for numerical inversion of Laplace transforms. The MATLAB code used in simulations of this paper, *invlap.m*, can freely be downloaded from http://www.mathworks.com/matlabcentral/fileexchange/ .

**Example 1.** The following transfer function appears in the one-link flexible robot arm [31]:

$$G(s) = \frac{-4.906s^2 - 0.5884s + 335.17}{s(s^3 + 0.55437s^2 + 139.6s + 27.91)}$$

$$= \frac{335.17\left(1 - \frac{s}{8.2057}\right)\left(1 + \frac{s}{8.3257}\right)}{s^4 + 0.55437s^3 + 139.6s^2 + 27.91s}. \quad (7)$$

The above transfer function has a non-minimum phase zero located at $z = 8.2057$ and four poles located at $p_1 = 0$, $p_2 = -0.2$, $p_{3,4} = -0.1772 \pm j11.8109$. Note that this system constitutes a relatively difficult control problem since it has a non-minimum phase zero and two complex-conjugate poles with a very low damping ratio ($\zeta = 0.0150$). Assuming

$$C_1(s) = \frac{1}{\sum_{k=1}^{v}\left(\frac{s}{8.2057}\right)^{\frac{k-1}{v}}}, \quad (8)$$

Yields

$$P(s) = \frac{335.17\left(1 - \left(\frac{s}{8.2057}\right)^{1/v}\right)\left(1 + \frac{s}{8.3257}\right)}{s^4 + 0.55437s^3 + 139.6s^2 + 27.91s}. \quad (9)$$

Since $P(s)$ has a pole at the origin, tracking of step command without steady-state error can be achieved simply by using a PD-type controller. In order to design the PD controller, $C_2(s)$, first we assign a value to $v$ and then tune the parameters of the controller assuming that the transfer function of process is equal to $P(s)$. Assuming $v = 20$, after a simple trial and error the transfer function of controller is obtained as the following:

$$C_2(s) = 0.1 + 0.5s. \quad (10)$$

(Note that the low-pass filter of derivative term is neglected for the sake of simplicity.) Fig. 3 shows the unit step response of the corresponding closed-loop system for $v = 15, 20, 25$. The very important observation in this figure is that the step response does not exhibit a sensible initial undershoot. In fact, since $G(s)$ (as well as the closed-loop transfer function) has odd number of non-minimum phase zeros, it is expected that mere application of any PID controller leads to initial undershoot in the step response. Hence, it can be concluded that using the fractional-order pole-zero canceller has the important property of decreasing the initial undershoots. A relevant discussion can be found in [18]. Note that in this example the final controller (using the nominal value of $v = 20$) is equal to the series connection of $C_1(s)$ and $C_2(s)$ as the following:

$$C(s) = C_1(s)C_2(s) = \frac{0.1 + 0.5s}{\sum_{k=1}^{20}\left(\frac{s}{8.2057}\right)^{(k-1)/20}}, \quad (11)$$

Which is almost bi-proper (the degree of numerator and denominator is equal to unity and 19/20, respectively).

Here, it should be emphasized that in general it is not necessary to use large values of $v$. In fact, in many cases even small values of $v$ lead to satisfactory results. For example, Fig. 4 shows the unit step response of the closed-loop system for $v = 2$ and $C_2(s) = 0.05 + 0.05s$. As it is observed in this figure, the response is satisfactory and still does not exhibit a sensible initial undershoot. However, it should be remind that increasing $v$ commonly increases the control effort.

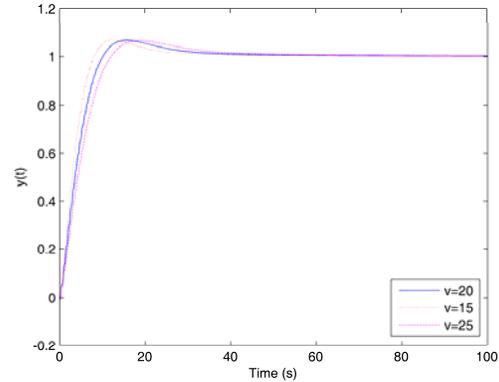

Fig. 3 Unit step response of the closed-loop system shown in Fig. 2 for different values of $v$ when the PD controller given in (10) is applied

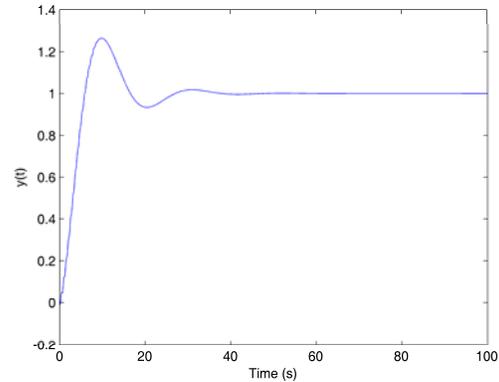

Fig. 4 Unit step response of the closed-loop system for $v$=2 and $C_2(s)$=0.05+0.05$s$, corresponding to Example 1

Finally, it should be noticed that neither the step response of Fig. 3 nor of Fig. 4 are obtained using optimal controllers and better responses can be obtained in both cases.

**Example 2.** The following transfer function is obtained by identification of a flexible-link manipulator [32]:

$$G(s) = \frac{b_6 s^6 + \ldots + b_1 s + b_0}{a_9 s^9 + \ldots + a_1 s + a_0}, \quad (12)$$

Where $a_9 = 1$, $a_8 = 486.7$, $a_7 = 69317.7$, $a_6 = 0.1616 \times 10^8$, $a_5 = 0.1062 \times 10^{10}$, $a_4 = 0.6167 \times 10^{11}$, $a_3 = 0.2624 \times 10^{13}$, $a_2 = 0.3595 \times 10^{14}$, $a_1 = 0.142 \times 10^{15}$, $a_0 = 0$, $b_6 = -14340.4953$, $b_5 = 0.4446 \times 10^7$, $b_4 = 0.5697 \times 10^9$, $b_3 = -0.1908 \times 10^{11}$, $b_2 = -0.9354 \times 10^{12}$, $b_1 = 0.6919 \times 10^{13}$, $b_0 = 0.2839 \times 10^{15}$. This system has three non-minimum phase zeros located at $z_1 = 400.0282$, $z_2 = 45.0015$, and $z_3 = 19.9982$. Moreover, similar to the previous example, it has poles with very low damping ratios. Since the system itself has a pole at the origin we design a PD-type controller. In this example we subject all of the non-minimum phase zeros of $G(s)$ to the fractional-order pole-zero cancellation assuming $v = 5$. That is, we consider $C_1(s)$ as the following:

$$C_1(s) = \frac{1}{\prod_{i=1}^{3} Q_{z_i,5}(s)}, \quad (13)$$

Where

$$Q_{z_i,5} = \sum_{k=1}^{5} \left(\frac{s}{z_i}\right)^{(k-1)/5}. \quad (14)$$

Then after a simple trial and error the corresponding PD controller is obtained as $C_2(s) = 5 + 2s$. (Note that similar to the previous example, this controller is not optimal in any sense and many other controllers can be designed instead. However, it is sufficient for the purpose of this example.) Fig. 5 shows the unit step response of the corresponding closed-loop system for $v = 4, 5, 6$.

Few points should be mentioned here. First, as it is observed in Fig. 5, the closed-loop system becomes faster (of course, at the cost of increasing undershoots and using a more control effort) by decreasing the value of $v$. It is a general observation that can be explained based on the relation between $v$ and bandwidth of the closed-loop system. Second, although $G(s)$ has an odd number of non-minimum phase zeros, no considerable undershoot is observed in the closed-loop step response, which is because of the effect of pre-compensator. The small undershoots observed in Fig. 5 can be decreased by changing the value of $v$ and parameters of the controller. The last point is that it is not necessary to cancel all of the non-minimum phase zeros of $G(s)$ to the same degree (here $v = 5$). In fact the performance of the closed-loop system can be better adjusted by suitable choice of $v_1$, $v_2$ and $v_3$.

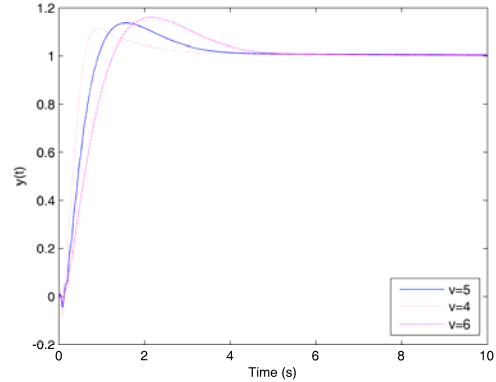

Fig. 5 Unit step response of the closed-loop system shown in Fig. 2 for different values of $v$ when the PD controller $C_2(s)=5+2s$ is applied